\documentclass[12pt]{article}
\usepackage{amsmath}
\usepackage{amssymb}
\usepackage{vatola}

\def\q{\quad}
\def\qq{\qquad}
\def\qtq#1{\q\t{#1}\q}
\def\[{\left[}
\def\]{\right]}
\def\t{\text}
\def\mod#1{\ (\text{\rm mod}\ #1)}
\def\f{\frac}
\def\e{\equiv}
\def\b{\binom}
\def\ex{\t{ex}}
\def\Ex{\t{Ex}}
\def\sls#1#2{(\f{#1}{#2})}

\let \pro=\proclaim
\let \endpro=\endproclaim

\begin{document}
\leftline{Colloquium Mathematicum 139(2015), no.2,
273-298}\par\q\newline
 \centerline {\bf Tur\'an's problem and Ramsey numbers for trees}
\par\q\newline \centerline{Zhi-Hong Sun$^1$, Lin-Lin Wang$^2$ and Yi-Li
Wu$^3$} \par\q\newline \centerline{$\ ^1$School of Mathematical
Sciences, Huaiyin Normal University,} \centerline{Huaian, Jiangsu
223001, P.R. China} \centerline{Email: zhihongsun@yahoo.com}
\centerline{Homepage: http://www.hytc.edu.cn/xsjl/szh}

\centerline{$\ ^2$School of Science, China University of Mining and
Technology,}
 \centerline{Xuzhou, Jiangsu  221116,
 PR China}
\centerline{Email: wanglinlin$\overline{\  }$1986@aliyun.com}

\centerline{$\ ^3$School of Mathematical Sciences, Huaiyin Normal
University,} \centerline {Huaian, Jiangsu 223001, PR China}
\centerline{Email: yiliwu@126.com}

\abstract{Let $T_n^1=(V,E_1)$ and $T_n^2=(V,E_2)$ be the trees on
$n$ vertices with
 $V=\{v_0,v_1,\ldots,v_{n-1}\}$,
 $E_1=\{v_0v_1,\ldots,v_0v_{n-3},v_{n-4}v_{n-2},v_{n-3}v_{n-1}\}$,
and $E_2=\{v_0v_1,\ldots,$ $v_0v_{n-3},v_{n-3}v_{n-2}, v_{n-3}$
$v_{n-1}\}$. In this paper, for $p\ge n\ge 5$ we obtain explicit
formulas for $\ex(p;T_n^1)$ and $\ex(p;T_n^2)$, where $\ex(p;L)$
denotes the maximal number of edges in a graph of order $p$ not
containing $L$ as a subgraph. Let $r(G\sb 1, G\sb 2)$ be the Ramsey
number of the two graphs $G_1$ and $G_2$. In this paper we also
obtain some explicit formulas for $r(T_m,T_n^i)$, where
$i\in\{1,2\}$ and $T_m$ is a tree on $m$ vertices with
$\Delta(T_m)\le m-3$.
\par\q
\newline 2010 Mathematics Subject Classification: 05C55, 05C35, 05C05
\newline Key words and phrases: Ramsey number, tree, Tur\'an's problem}
\endabstract

\section*{1. Introduction}
 \par\q\  In this paper, all graphs are simple
graphs. For a graph $G=(V(G),E(G))$ let $e(G)=|E(G)|$ be the number
of edges in $G$ and let $\Delta(G)$ be the maximal degree of $G$.
 For a forbidden graph $L$, let $\ex(p;L)$
denote the maximal number of edges in a graph of order $p$ not
containing any copies of $L$. The corresponding Tur\'an problem is
to evaluate $\ex(p;L)$. For a graph $G$ of order $p$, if $G$ does
not contain any copies of $L$ and $e(G)=\ex(p;L)$, we say that $G$
is an extremal graph. In this paper we also use $\Ex(p;L)$ to denote
the set of extremal graphs of order $p$ not containing $L$ as a
subgraph.
\par Let $\Bbb N$ be the set of
positive integers. Let $p,n\in\Bbb N$ with $p\ge n\ge 2$. For a
given tree $T_n$ on $n$ vertices,
 it is difficult to determine the value of $\ex(p;T_n)$.
The famous Erd\H os-S\'os
 conjecture asserts that
$\ex(p;T_n)\le \f{(n-2)p}2$. For the progress on the Erd\H os-S\'os
 conjecture, see for example [8, 11].
 Write
$p=k(n-1)+r$, where $k\in\Bbb N$ and $r\in\{0,1,\ldots,n-2\}$. Let
$P_n$ be the path on $n$ vertices. In [4] Faudree and Schelp showed
that
$$\ex(p;P_n)=k\binom {n-1}2+\binom r2=\f{(n-2)p-r(n-1-r)}2.\tag 1.1$$
 Let $K_{1,n-1}$ denote the unique tree on $n$
vertices with $\Delta(K_{1,n-1})=n-1$, and let
  $T_n'$ denote the unique tree on $n$ vertices
with $\Delta(T_n')=n-2$. For $n\ge 4$ let $T_n^*=(V,E)$ be the tree
on $n$ vertices with
 $V=\{v_0,v_1,\ldots,v_{n-1}\}$ and
 $E=\{v_0v_1,\ldots,v_0v_{n-3},v_{n-3}v_{n-2},$ $v_{n-2}v_{n-1}\}$.
 In [10] we determine $\ex(p;K_{1,n-1})$, $\ex(p;T_n')$
and $\ex(p;T_n^*)$. For $i=1,2$ let $T_n^i=(V,E_i)$ be the tree on
$n$ vertices with
 $$\align &V=\{v_0,v_1,\ldots,v_{n-1}\},
 \\&E_1=\{v_0v_1,\ldots,v_0v_{n-3},v_{n-4}v_{n-2},v_{n-3}v_{n-1}\},
\\&E_2=\{v_0v_1,\ldots,v_0v_{n-3},v_{n-3}v_{n-2},v_{n-3}v_{n-1}\}.
\endalign$$
In this paper, for $p\ge n\ge 5$ we obtain explicit formulas for
$\ex(p;T_n^1)$ and $\ex(p;T_n^2)$ (see Theorems 2.1 and 3.1).

\par For a graph $G$, as usual $\overline{G}$ denotes the
complement of $G$. Let $G\sb 1$ and $G\sb 2$ be two graphs. The
Ramsey number $r(G\sb 1, G\sb 2)$ is the smallest positive integer
$p$ such that, for every graph $G$ with $p$ vertices, either $G$
contains  a copy of $G\sb 1$
 or else $\overline{G}$ contains a copy of $G_2$.
 \par Let $n\in\Bbb N$, $n\ge 6$ and let $T_n$ be a tree on $n$ vertices. As mentioned in [7],
 recently Zhao proved the following conjecture of
 Burr and Erd\H os [2]:  $r(T_n,T_n)\le 2n-2$.
 Let $m,n\in\Bbb N$. In 1973 Burr and
 Roberts [3] showed that for $m,n\ge 3$,
$$r(K_{1,m-1},K_{1,n-1})=\cases m+n-3&\t{if $2\nmid mn$,}
\\m+n-2&\t{if $2\mid mn$}.\endcases\tag 1.2$$
In 1995, Guo and Volkmann [5] proved that for $n>m\ge 4$,
$$r(K_{1,m-1},T_n')=\cases m+n-3&\t{if $2\mid m(n-1)$,}
\\m+n-4&\t{if $2\nmid m(n-1)$.}
\endcases\tag 1.3$$
Recently the first author evaluated the Ramsey number $r(T_m,T_n^*)$
for  $T_m\in\{P_m,K_{1,m-1},$ $T_m',T_m^*\}$. In particular, he
proved that (see [9]) for $n>m\ge 7$,
    $$r(K_{1,m-1},T_n^*)=\cases m+n-3&\t{if $m-1\mid n-3$,}
\\m+n-4&\t{if $m-1\nmid n-3$.}\endcases\tag 1.4$$

\par  Suppose
$m,n\in\Bbb N$ and $i,j\in\{1,2\}$. In this paper, using the formula
for $\ex(p;T_n^i)$ and the method in [9] we evaluate $r(T_m,T_n^i)$
for $T_m\in\{K_{1,m-1}, T_m',T_m^*,$ $T_m^j\}$. In particular, we
have the following typical results:
$$\aligned
&r(T_n^i,T_n^j)=2n-6-(1-(-1)^n)/2,\ r(P_n,T_n^j)=2n-7\qtq{for}n\ge
17,
\\&r(T_n^i,T_n')=r(T_n^i,T_n^*)=2n-5\qtq{for}n\ge 8,
\\&r(K_{1,m-1},T_n^i)=m+n-4\qtq{for} n>m\ge 7
\q \t{and}\q 2\mid mn,
\\&r(T_m^i,T_n^j)=m+n-5\qtq{for}m\ge 7,\ n\ge
(m-3)^2+3\qtq{and}m-1\nmid n-4,\endaligned$$ and for $n>m\ge 16$,
$$ r(T_m',T_n^i)=\cases m+n-4&\t{if $m-1\mid n-4$,}
\\m+n-6&\t{if $n=m+1\e 1\mod 2$,}
\\m+n-5&\t{otherwise.}\endcases$$
\par
In addition to the notation introduced above, throughout the paper
we also use the following symbols: $[x]$ is the greatest integer not
exceeding $x$, $d(v)$ is the degree of the vertex $v$ in a graph,
$\Gamma(v)$ is the set of vertices adjacent to the vertex $v$,
$d(u,v)$ is the distance between the two vertices $u$ and $v$ in a
graph, $K_n$ is the complete graph on $n$ vertices, $G[V_0]$ is the
subgraph of $G$ induced by vertices in the set $V_0$ (we write
$G[v_1,\ldots,v_m]$ instead of $G[\{v_1,\ldots,v_m\}]$), $G-V_0$ is
the subgraph of $G$ obtained by deleting vertices in $V_0$ and all
edges incident to them, and finally $e(V_1V_1')$ is the number of
edges with one endpoint in $V_1$ and another endpoint in $V_1'$.
 \section*{2.
Evaluation of $\ex(p;T_n^1)$}
 \pro{\q\ Lemma 2.1} Let $p,n\in\Bbb N$
with $p\geq n-1\geq 1$. Then $\ex(p;K_{1,n-1})=[\f{(n-2)p}2]$.
\endpro
\par This is a known result. See for example [10, Theorem 2.1].
 \pro{Lemma 2.2} Let
$p,n\in \Bbb N$, $p\ge n\ge 7$ and  $G\in \Ex(p;T_n^1)$. Suppose
that $G$ is connected. Then $\Delta(G)=n-4$ and
$e(G)=[\f{(n-4)p}2]$.
\endpro Proof.
Since a graph not containing $K_{1,n-3}$ as a subgraph implies that
the graph does not contain $T_n^1$ as a subgraph, by Lemma 2.1 we
have
$$e(G)=\ex(p;T_n^1)\ge \ex(p;K_{1,n-3})=\Big[\f{(n-4)p}2\Big].\tag 2.1$$
If $\Delta(G)\le n-5$, using Euler's theorem we see that $e(G)=\f
12\sum_{v\in  V(G)}d(v)\le \f{(n-5)p}2.$ This together with (2.1)
yields $\f{(n-4)p-1}2\le[\f{(n-4)p}2]\le e(G)\le \f{(n-5)p}2$. This
is impossible. Hence $\Delta(G)\ge n-4$. Now we show that
$\Delta(G)=n-4$.
\par Suppose $q\ge n$ and $q=k(n-1)+r$ with
$k\in\Bbb N$ and
 $r\in\{0,1,\ldots,n-2\}$. Then clearly $kK_{n-1}\cup K_r$ does not
 contain any copies of $T_n^1$ and so $\ex(q;T_n^1) \ge
 e(kK_{n-1}\cup K_r)$.
 For $q=n$ we see that $e(kK_{n-1}\cup K_r)=e(K_{n-1}\cup K_1)=\f{(n-1)(n-2)}2>2n-1$.
 For $q\ge n+1$ we see that $(n-6)q\ge (n-6)(n+1)>\sls{n-1}2^2-2$
 and so $e(kK_{n-1}\cup K_r)=\f{k(n-1)(n-2)}2
 +\f{r(r-1)}2=\f{(n-2)q-r(n-1-r)}2\ge \f{(n-2)q-(\f{n-1}2)^2}2>2q-1$.
 Hence
$$\ex(q;T_n^1)\ge e(kK_{n-1}\cup K_r)>2q-1\qtq{for}q\ge n.\tag 2.2$$
\par
Suppose $v_0\in V(G), d(v_0)=\Delta(G)=m$ and
$\Gamma(v_0)=\{v_1,\ldots,v_m\}$. If  $m=p-1$, as $G$ does not
contain $T_n^1$ as a subgraph, we see that $G[v_1,\ldots,v_m]$ does
not contain $2K_2$ as a subgraph and hence $e(G[v_1,\ldots,v_m])\le
m-1$. Therefore
 $$e(G)=d(v_0)+e(G[v_1,\ldots,v_m])\le m+m-1=2p-3.\tag 2.3$$
  By (2.2), we have $e(G)=\ex(p;T_n^1)>2p-1$ and we get
   a contradiction. Hence $m<p-1$. Suppose that
$u_1,\ldots,u_t$ are all vertices in $G$ such that
$d(u_1,v_0)=\cdots=d(u_t,v_0)=2$. Then $t\ge 1$. Assume $u_1v_1\in
E(G)$ with no loss of generality.  If $m=p-2,$ then
$V(G)=\{v_0,v_1,\ldots, v_m,u_1\}$ and $v_iv_j\not\in E(G)$ for
$2\le i<j\le m$. If $v_1v_i\in E(G)$ for some
$i\in\{2,3,\ldots,m\}$, then $u_1v_j\not\in E(G)$ for all
$j\not=1,i$. Hence $\ex(p;T_n^1)=e(G)\le \max\{2m,m+3\}\le 2m=2p-4$,
which contradicts (2.2).
\par By the above, $m<p-2$. We first assume $m\ge n-2$. As $G$ does not contain any copies of $T_n^1$, we see
that $\{v_2,\ldots,v_m\}$ is an independent set, $u_iv_j\not\in
E(G)$ for any $i\in\{2,3,\ldots,t\}$ and $j\in\{2,3,\ldots,m\}$, and
$u_iv_1\in E(G)$ for any $i=1,2,\ldots,t$.
 Set $V_1=\{v_0,v_2,v_3,\ldots,v_m\}$. Then $e(G[V_1])=m-1$.
If $u_1$ is adjacent to at least two vertices in
$\{v_2,v_3,\ldots,v_m\}$, then $v_1v_j\notin E(G)$ for any
$j=2,3,\ldots,m$. If $v_1$ is adjacent to at least two vertices in
$\{v_2,v_3,\ldots,v_m\}$, then $u_1v_j\notin E(G)$ for any
$j=2,3,\ldots,m$.
 Hence there are at most $m$ edges with one endpoint in $V_1$ and
 another endpoint in $G-V_1$. Therefore,
$$e(G)\le e(G[V_1])+m+e(G-V_1)=2m-1+e(G-V_1).\tag 2.4$$
For $m\in\{n-2,n-1\}$ let $G_1=K_m$. Then clearly
$e(G_1)=\f{m(m-1)}2> 2m-1$. For $m=k(n-1)+r\ge n$ with $k\in\Bbb N$
and $0\le r\le n-2$ let $G_1=kK_{n-1}\cup K_r$. Then $G_1$ does not
contain any copies of $T_n^1$ and $e(G_1)>2m-1$ by (2.2). Thus, by
(2.4) we have $e(G)\le 2m-1+e(G-V_1)<e(G_1\cup (G-V_1))$ for $m\ge
n-2$. This contradicts the fact $G\in \Ex(p;T_n^1)$.
\par Suppose $m=n-3$ and $d(v_1)=n-3$.
Then $v_1v_s\not\in E(G)$ for some $s\in\{2,3,\ldots,n-3\}$. We
claim that $V(G)=\{v_0,v_1,\ldots,v_m,u_1, \ldots,$ $u_t\}.$
Otherwise, there exists  $w\in V(G)$ such that $d(v_0,w)=3.$ As
$d(v_1)=n-3$, we see that the subgraph induced by $\{v_1,v_s,w\}\cup
\Gamma(v_1)$ contains a copy of $T_n^1$. This contradicts the
assumption $G\in \Ex(p;T_n^1)$. Hence the claim is true and so
$|V(G)|=p=n-2+t$. Since $p\ge n$ we have $t\ge 2$. For
$i=1,2,\ldots,t$ and $j=2,3,\ldots,n-3$ we have $u_iv_j\not\in
E(G)$, $u_iv_1\in E(G)$ and so $t+1\le d(v_1)=n-3.$ Therefore $2\le
t\le n-4$ and hence
$$\align e(G)&=e(G[v_0,v_2,v_3,\ldots,v_{n-3}])+d(v_1)+e(G[u_1,\ldots,u_t])
\\&\le\b{n-3}2+n-3+\b t2=\b{n-2}2+\b t2.\endalign$$
Clearly $K_{n-1}\cup K_{t-1}$ does not contain $T_n^1$ and
$$e(K_{n-1}\cup K_{t-1})=\b{n-1}2+\b{t-1}2=\b{n-2}2+\b t2+n-1-t>e(G).$$
This contradicts the assumption $G\in \Ex(n-2+t;T_n^1)$.
\par Now suppose $m=n-3$ and $d(v_1)\le n-4$.
If $t=1$, setting $V_2=\{v_0,v_1,\ldots,v_{n-3},u_1\}$ we see that
$$\align
e(G)&=e(G[v_0,v_2,v_3,\ldots,v_{n-3}])+d(v_1)+d(u_1)-1+e(G-V_2)
\\&\le \b{n-3}2+n-4+n-4+e(G-V_2)
\\&=\f{n^2-3n-4}2+e(G-V_2)<e(K_{n-1}\cup (G-V_2)).\endalign$$
This contradicts the assumption $G\in \Ex(p;T_n^1)$. Hence $t\ge 2.$
For $i=1,2,\ldots,t$ and $j=2,3,\ldots,n-3$ we see that
$u_iv_j\not\in E(G)$ and $u_iv_1\in E(G)$. Let
$V_3=\{v_0,v_1,\ldots,v_{n-3}\}.$ Then $$\align
e(G)&=d(v_1)+e(G[v_0,v_2,v_3,\ldots,v_{n-3}])+e(G-V_3)
\\&\le n-4+\b{n-3}2+e(G-V_3)=\f{n^2-5n+4}2+e(G-V_3)
\\&<e(K_{n-2}\cup (G-V_3)).\endalign$$ Since $G$ is an extremal graph,
we get  a contradiction.
\par Summarizing all the above we obtain $\Delta(G)=n-4$
and so $e(G)=\sum_{v\in V(G)}d(v)\le \f{(n-4)p}2$. This together
with (2.1) yields $e(G)=[\f{(n-4)p}2]$, which completes the proof.

\pro{Lemma 2.3} Let $n,n_1,n_2\in\Bbb N$ with $n_1<n-1$ and
$n_2<n-1$. Then
$$\b{n_1}2+\b{n_2}2<\min\Big\{\b{n_1+n_2}2,
\b{n-1}2+\b{n_1+n_2-n+1}2\Big\}.$$
\endpro
Proof. It is clear that
$$\b{n_1}2+\b{n_2}2=\f{(n_1+n_2)(n_1+n_2-1)
-2n_1n_2}2<\b{n_1+n_2}2$$ and
$$\align &\b{n-1}2+\b
{n_1+n_2-n+1}2-\b{n_1}2-\b{n_2}2
\\&=\f{(n-1)(n-2)+(n_1+n_2-n+1)(n_1+n_2-n)}2-\f{(n_1+n_2)(n_1+n_2-1)-2n_1n_2}2
\\&=(n-1-n_1)(n-1-n_2)>0.\endalign$$
Thus the lemma is proved.

\pro{Lemma 2.4} Suppose that $p\in\Bbb N$, $p\ge 6$, and G is a
connected graph of order p that does not contain any copies of
$T_6^1$. Then $e(G)\le{2p-3}$.\endpro
 Proof. Clearly $\Delta(T_6^1)=3$. Suppose $v_0\in V(G)$,
$d(v_0)=\Delta(G)=m$ and $\Gamma(v_0)=\{v_1,\ldots,v_m\}$. If
$\Delta(G)=m\le3$, using Euler's theorem we see that
$e(G)\le\frac{3p}2\le2p-3$. From now on we assume
$\Delta(G)=m\geq4$. If $d(v)\le 2$ for all $v\in V(G)-\{v_0\}$, then
$$ e(G)=\frac 12\sum_{v\in V(G)}d(v)\le\frac
12\big(m+2(p-1)\big)\le \f{3(p-1)}2<2p-3.$$  So the result is true.
Now we assume $d(v)\ge 3$ for some $v\in V(G)-\{v_0\}$.  We may
choose a vertex $u_0\in V(G)$ so that $u_0\neq v_0$, $d(u_0)\ge 3$
and $d(u_0,v_0)$ is as small as possible.

\par We first assume $d(u_0,v_0)=1$ and $u_0=v_1$ with no loss of
generality. That is, $d(v_1)\ge 3$. Suppose $\Gamma(v_1)\subset
\{v_0,v_1,\ldots,v_m\}$. Since $d(v_1)\ge 3$ and $G$ does not
contain any copies of $T_6^1$, we see that
$V(G)=\{v_0,\ldots,v_m\}$, $m=p-1\ge 5$ and $G[v_1,\ldots,v_m]$ does
not contain any copies of $2K_2$. Thus $e(G)\le d(v_0)+m-1=2m-1\le
2(m+1)-3=2p-3$. Now assume
$\Gamma(v_1)-\{v_0,v_1,\ldots,v_m\}=\{w_1,\ldots,w_t\}$. Since
$d(v_0)=m\ge 5$, $d(v_1)\ge 3$ and $G$ does not contain any copies
of $T_6^1$, we see that $V(G)=\{v_0,v_1,\ldots,v_m,w_1,\ldots,w_t\}$
and $\{v_2,\ldots,v_m\}$ is an independent set.
 For
$t\ge 2$, we have $e(G[w_1,\ldots,w_t])\le 1$ and $v_iw_j\notin
E(G)$ for any $i \in\{2,3,\ldots,m\}$ and $j \in\{1,2,\ldots,t\}$.
Thus $e(G)\le d(v_0)+d(v_1)-1+1\le 2m<2(m+1+t)-3=2p-3.$ Now assume
$t=1$. Then $v_1v_i\in E(G)$ for some $i\in\{2,3,\ldots,m\}$ and
$v_jw_1\notin E(G)$ for $j\in\Bbb\{2,3,\ldots,m\}-\{i\}$. Hence
$e(G)\le d(v_0)+d(v_1)-1+1\le 2m<2(m+2)-3=2p-3$.

\par Next we assume
$d(u_0,v_0)=2$. Then  $\{v_1,\ldots,v_m\}$ is an independent set. If
$\Gamma(u_0)\subseteq\{v_1,\ldots,v_m\}$, then
$V(G)=\{v_0,\ldots,v_m,u_0\}$ and so $e(G)=d(v_0)+d(u_0)\le
m+m<2(m+2)-3=2p-3$. If
$\Gamma(u_0)-\{v_2,\ldots,v_m\}=\{v_1,w_1,\ldots,w_t\}$, we see that
$V(G)=\{v_0,v_1,\ldots,v_m,u_0,$ $w_1,\ldots,w_t\}$ and so $e(G)
=d(v_0)+d(u_0)+e(G[w_1,\ldots,w_t])\le m+m+1<2(m+2+t)-3=2p-3$.

\par Finally we assume $d(u_0,v_0)\ge 3$. Suppose that
$v_0v_1u_1u_2\cdots u_ku_0$ is the shortest path in $G$ between
$v_0$ and $u_0$, and $\Gamma(u_0)=\{w_1,\ldots,w_t,u_k\}$. Since $G$
is connected and $G$ does not contain any copies of $T_6^1$, it is
easily seen that $V(G)=\{v_0,v_1,\ldots,v_m,u_1,\ldots,u_k,u_0,$
$w_1,\ldots,w_t\}$, $d(v_2)=\cdots=d(v_m)=1$, $d(v_1)=d(u_1)=\cdots=
d(u_k)=2$ and
 $e(G[w_1,\ldots,w_t])\le 1$.
  Clearly $G$ is a tree
or a graph obtained by adding an edge to a tree. Hence $e(G)\le
p<2p-3$.
\par Summarizing all the above proves the lemma.

\pro{Theorem 2.1} Suppose $p,n\in\Bbb N,\ p\geq n-1\geq 4$ and
$p=k(n-1)+r$, where $k\in \Bbb N$ and $r\in\{ 0,1,\ldots,n-2\}$.
Then $$\aligned \ex(p;T_n^1)&=
\max\Big\{\Big[\f{(n-2)p}2\Big]-(n-1+r),\f{(n-2)p-r(n-1-r)}2\Big\}
\\&= \cases [\f{(n-2)p}2]-(n-1+r) &\t{if}\ n\ge 16\ \t{and}\
3\le r\le n-6 \ \t{or if}\\&\ \ 13\le n\le 15 \ \t{and}\ 4\le r\le n-7,\\
\f{(n-2)p-r(n-1-r)}2&\t{otherwise.}\endcases\endaligned$$
\endpro
Proof. Clearly $\ex(n-1;T_n^1)=e(K_{n-1})=\f{(n-2)(n-1)}2$. Thus the
result is true for $p=n-1$. From now on we assume $p\ge n$. Since
$T_5^1\cong P_5$, by (1.1) we obtain the result in the case $n=5$.
Now we assume $n\ge 6$. Suppose $G\in \Ex(p;T_n^1)$ and
$G_1,\ldots,G_t$ are all components of $G$ with $|V(G_i)|=p_i$ and
$p_1\le p_2\le\cdots\le p_t.$ Then clearly $G_i\in \Ex(p_i;T_n^1)$
for $i=1,2,\ldots,t$.
\par We first consider the case $n=6$. If $p_i\le 5$, then clearly
$G_i\cong K_{p_i}$ and $e(G_i)=\b{p_i}2$. If $p_i\ge 6$ and
$p_i=5k_i+r_i$ with $k_i\in\Bbb N$ and $0\le r_i\le 4$, from Lemma
2.4 we have $e(G_i)\le 2p_i-3\le 2p_i-\f{r_i(5-r_i)}2=e(k_iK_5\cup
K_{r_i}).$ Since $k_iK_5\cup K_{r_i}$ does not contain any copies of
$T_6^1$ and $G_i\in \Ex(p_i;T_6^1)$, we see that $e(G_i)\ge
e(k_iK_5\cup K_{r_i})$ and so $e(G_i)=e(k_iK_5\cup K_{r_i})$.
Therefore, there is a graph $G'\in \Ex(p;T_6^1)$ such that
$G'=a_1K_1\cup a_2K_2\cup a_3K_3\cup a_4K_4\cup a_5K_5$, where
$a_1,\ldots, a_5$ are nonnegative integers. If $a_1+a_2+a_3+a_4\le
1$, then $\ex(p;T_6^1)=e(G')=e(a_5K_5\cup K_r)=k\b 52+\b r2$. If
$a_1+a_2+a_3+a_4>1$, then $2a_1+3a_2+3a_3+2a_4>3\ge \f{r(5-r)}2$ and
so
$$\align
&e(a_1K_1\cup a_2K_2\cup a_3K_3\cup a_4K_4)
\\&=a_2+3a_3+6a_4<2(a_1+2a_2+3a_3+4a_4)-\f{r(5-r)}2=(k-a_5)\b
52+\b r2.\endalign$$ Thus, $\ex(p;T_6^1)= e(G')=e(a_1K_1\cup
a_2K_2\cup a_3K_3\cup a_4K_4)+e(a_5K_5)<k\b 52+\b r2$. Since
$kK_5\cup K_r$ does not contain any copies of $T_6^1$, we get a
contradiction. Thus $\ex(p;T_6^1)=e(kK_5\cup K_r)=k\b 52+\b
r2=2p-\f{r(5-r)}2$. This proves the result for $n=6$.
\par From now on we assume
$n\ge 7$. If $t=1$, then $G$ is connected. Thus,
  by Lemma 2.2
we have
 $$e(G)=\Big[\f{(n-4)p}2\Big]\qtq{for}t=1.\tag 2.5 $$
   \par Now we assume $t\ge 2$. We claim
that $p_i\ge n-1$ for $i\ge 2$. Otherwise, $p_1\le p_2<n-1$ and so
$G_1\cup G_2\cong K_{p_1}\cup K_{p_2}.$ If $p_1+p_2<n,$ by Lemma 2.3
we have
 $e(G_1\cup G_2)=e(K_{p_1}\cup
 K_{p_2})=\b{p_1}2+\b{p_2}2<\b{p_1+p_2}2=e(K_{p_1+p_2}).$
Since $K_{p_1+p_2}$ does not contain $T_n^1$ and $G_1\cup G_2\in
\Ex(p_1+p_2;T_n^1)$ we get a contradiction. Hence $p_1+p_2\ge n$.
Using Lemma 2.3 again we see that
 $$\align
 e(G_1\cup G_2)&=e(K_{p_1}\cup K_{p_2})=\b{p_1}2+\b{p_2}2
 \\&<\b{n-1}2+\b{p_1+p_2-n+1}2=e(K_{n-1}\cup K_{p_1+p_2-n+1}).
 \endalign$$
Since $p_1\le p_2<n-1,$ we have $p_1+p_2-n+1<n-1$. Therefore
$K_{n-1}\cup K_{p_1+p_2-n+1}$ does not contain $T_n^1$. As $G_1\cup
G_2$ is an extremal graph without $T_n^1$, we also get a
contradiction. Thus, the claim is true.
 \par Next we claim that
$p_i\le n-1$ for all $i=1,2,\ldots,t-1$. If $p_{t-1}\ge n$, by Lemma
2.2 we have
$$e(G_{t-1}\cup G_t)=e(G_{t-1})+e(G_t)
=\[\f{(n-4)p_{t-1}}2\]+\[\f{(n-4)p_t}2\] \le
\[\f{(n-4)(p_{t-1}+p_t)}2\].$$
Let $H\in \Ex(p_{t-1}+p_t-n+1;K_{1,n-3})$. As $p_{t-1}+p_t-n+1\ge
p_t+1\ge n+1$, we have $e(H)=[\f{(n-4)(p_{t-1}+p_t-n+1)}2]$ by Lemma
2.1. Clearly $K_{n-1}\cup H$ does not contain any copies of $T_n^1$
and
$$\align e(K_{n-1}\cup H)&
=e(K_{n-1})+e(H)=\b{n-1}2+\[\f{(n-4)(p_{t-1}+p_t-n+1)}2\]
\\&=\[\f{(n-4)(p_{t-1}+p_t)}2\]+n-1> e(G_{t-1}\cup G_t).\endalign$$
Since $G_{t-1}\cup G_t\in \Ex(p_{t-1}+p_t;T_n^1)$, we get a
contradiction. Hence $p_1\le p_2\le\cdots \le p_{t-1}\le n-1$.
Combining this with the previous assertion that $p_t\ge \cdots \ge
p_2\ge n-1$ we obtain
 $$p_1\le n-1,\q p_2=\cdots=p_{t-1}=n-1\qtq{and}
p_t\ge n-1.\tag 2.6$$ As $G$ is an extremal graph, we must have
$$G_1\cong K_{p_1},\q G_2\cong K_{n-1},\q\ldots, \q G_{t-1}\cong
K_{n-1}.\tag 2.7$$

\par If $p_t=n-1$, then $G_t\cong K_{n-1}$.
 By (2.7), $G\cong K_{p_1}\cup (t-1)K_{n-1}\cong kK_{n-1}\cup K_r$.
Thus,
$$e(G)=k\b{n-1}2+\b r2
=\f{(n-2)p-r(n-1-r)}2\ \t{for}\ t\ge 2\ \t{and}\ p_t=n-1.\tag 2.8 $$
\par Now we assume  $p_t\ge n$. By Lemma 2.2,
$e(G_t)=[\f{(n-4)p_t}2]$. Since $p_1\le n-1$,
 we have $G_1\cong K_{p_1}$ and so $e(G_1)=e(K_{p_1})=\b{p_1}2.$ Let
$H_1\in \Ex(p_1+p_t;K_{1,n-3}).$ Then $H_1$ does not contain $T_n^1$
as a subgraph. By Lemma 2.1, for $p_1\le n-4$ we have $$\align
e(H_1)&=\[\f{(n-4)(p_1+p_t)}2\] \ge
\[\f{(n-4)p_t}2\]+\[\f{(n-4)p_1}2\]
\\&\ge  \[\f{(n-4)p_t}2\]+\f{(n-4)(p_1-1)}2+1
\\&> \[\f{(n-4)p_t}2\]+\f{p_1(p_1-1)}2
=e(G_1\cup G_t).\endalign$$ This contradicts $G_1\cup
 G_t\in \Ex(p_1+p_t;T_n^1)$.
Hence $n-3\le p_1\le n-1.$
\par For $p_1\in\{n-3,n-2\}$ and $p_t\ge n$, we have
$p_1(p_1-(n-3))\le 2n-4$ and so
$$\align e(G_1\cup G_t)&=e(G_1)+e(G_t)=\b{p_1}2+\[\f{(n-4)p_t}2\]
\\&\le \f{p_1(p_1-1)+(n-4)p_t}2=\f{p_1(p_1-(n-3))+(n-4)(p_1+p_t)}2
\\&\le \f{2n-4+(n-4)(p_1+p_t)}2=\b{n-1}2+\f{(n-4)(p_1+p_t-n+1)-2}2
\\&<\b{n-1}2+\Big[\f{(n-4)(p_1+p_t-n+1)}2\Big].\endalign$$
Let $H_2\in \Ex(p_1+p_t-n+1;K_{1,n-3})$. Then $K_{n-1}\cup H_2$ does
not contain any copies of $T_n^1$. Since $p_1+p_t-n+1\ge p_1+1\ge
n-2$, applying Lemma 2.1 we have $e(H_2)=[\f{(n-4)(p_1+p_t-n+1)}2]$.
 Thus, we have $e(K_{n-1}\cup H_2)=\b{n-1}2+[\f{(n-4)(p_1+p_t-n+1)}2]>e(G_1\cup
G_t)$. This contradicts $G_1\cup G_t\in \Ex(p_1+p_t;T_n^1)$.

\par By the above, for $t\ge 2$ and $p_t\ge n$ we have $p_1=p_2=\cdots
=p_{t-1}=n-1$. If $p_t\ge 2n-2$, setting $H_3\in
\Ex(p_t-(n-1);K_{1,n-3})$ and then applying Lemmas 2.1 and 2.2 we
find that
$$e(G_t)=\[\f{(n-4)p_t}2\]<\b{n-1}2+\[\f{(n-4)(p_t-(n-1))}2\]
=e(K_{n-1}\cup H_3).$$ This contradicts the fact $G_t\in
\Ex(p_t;T_n^1)$. Hence $n\le p_t<2n-2$ and so $r\ge 1$.
 Note that $p=k(n-1)+r=(k-1)(n-1)+n-1+r$ and $n\le n-1+r<2n-2$.
Hence $t=k$, $p_t=n-1+r$ and therefore
 $$\aligned e(G)&=e((k-1)K_{n-1})+e(G_t)
 =(k-1)\b{n-1}2+\[\f{(n-4)(n-1+r)}2\]
\\&= \[\f{(n-2)p}2\]-(n-1+r)\q \t{for $t\ge 2$ and $p_t\ge n$}.
\endaligned\tag 2.9$$
\par Since $G\in \Ex(p;T_n^1)$, by comparing (2.5), (2.8) and (2.9) we get
  $$e(G)=\max\bigg\{\[\f{(n-4)p}2\],
  \f{(n-2)p-r(n-1-r)}2,\[\f{(n-2)p}2\]-(n-1+r)\bigg\}.$$
Observe that $p=k(n-1)+r\ge n-1+r$. We see that
$[\f{(n-4)p}2]=[\f{(n-2)p}2]-p\le [\f{(n-2)p}2]-(n-1+r)$ and
therefore
$$\aligned \ex(p;T_n^1)&=e(G)=\max\bigg\{
\f{(n-2)p-r(n-1-r)}2,\[\f{(n-2)p}2\]-(n-1+r)\bigg\}
\\&=\f{(n-2)p-r(n-1-r)}2+\max\bigg\{0,\Big[\f{r(n-3-r)-2(n-1)}2\Big]\bigg\}.\endaligned\tag
2.10$$
 \par For $7\le n\le 12$ we have
 $r(n-3-r)-2(n-1)\le \f{(n-3)^2}4-2(n-1)=\f{(n-7)^2-32}4<0.$
For $r\in\{0,1,2,n-5,n-4,n-3,n-2\}$ we see that
 $r(n-3-r)-2(n-1)<0.$
Suppose $n\geq 13$ and $3\leq r\leq n-6.$ For $4\le r\le n-7$ we
have $|r-\f{n-3}2|\le \f{n-11}2$ and so
$$\align
r(n-3-r)-2(n-1)&=\f{n^2-14n+17}4-\Big(r-\f{n-3}2\Big)^2
\\&\geq \f{n^2-14n+17}4-\Big(\f{n-11}2\Big)^2=2n-26\ge
0.\endalign$$  For $r\in\{3,n-6\}$ we have
  $r(n-3-r)-2(n-1)=3(n-6)-2(n-1)=n-16.$
Now combining the above with (2.10) we deduce the result.
   \pro{Corollary 2.1}
Suppose $p,n\in\Bbb N$, $p\ge n\ge 5$ and $n-1\nmid p$. Then
$\f{(n-2)p}2-\f{(n-1)^2}8\le \ex(p;T_n^1)\le \f{(n-2)(p-1)}2$.
\endpro
Proof. Suppose $p=k(n-1)+r$ with $k\in\Bbb N$ and
$r\in\{0,1,\ldots,n-2\}$. Then $r\ge 1$. Clearly $\f{(n-1)^2}4\ge
r(n-1-r)=(\f{n-1}2)^2-(\f{n-1}2-r)^2\ge (\f{n-1}2)^2-(\f{n-1}2-1)^2
=n-2$ and $n-1+r>\f{n-2}2$. Thus, from Theorem 2.1 we deduce that
$\ex(p;T_n^1)\le \f{(n-2)p-(n-2)}2$ and $\ex(p;T_n^1)\ge
\f{(n-2)p-r(n-1-r)}2\ge \f{(n-2)p-(n-1)^2/4}2$. This proves the
corollary.
 \section*{3. Evaluation
of $\ex(p;T_n^2)$}
 \par\q\pro{Lemma 3.1} Let
$p,n\in\Bbb N$, $p\geq n\ge 7$ and $G\in \Ex(p;T_n^2)$. Suppose that
$G$ is connected. Then $\Delta(G)\le n-3$. Moreover, for $p<2n-2$ we
have $\Delta(G)\le n-4$.
 \endpro
  Proof. Since a graph does not contain $K_{1,n-3}$ implies that the graph
does not contain $T_n^{2}$, by Lemma 2.1 we have
$$e(G)=\ex(p;T_n^{2})\ge \ex(p;K_{1,n-3})=\Big[\f{(n-4)p}2\Big].
\tag 3.1$$
Suppose that $v_0\in V(G)$, $d(v_0)=\Delta(G)=m$ and
$\Gamma(v_0)=\{v_1,\ldots, v_m\}$. If $V(G)=\{v_0,v_1,\ldots,v_m\}$,
then $m=p-1\ge n-1$. Since $G$ does not contain $T_n^{2}$, we see
that $G[v_1,\ldots,v_m]$ does not contain $K_{1,2}$ and hence
$e(G[v_1,\ldots,v_m])\le \f m2$. Therefore
 $e(G)=d(v_0)+e(G[v_1,\ldots,v_m])\le m+\f m2=\f{3(p-1)}2\le
 \f{(n-4)p-3}2 <[\f{(n-4)p}2]. $
This contradicts (3.1). Thus $p>m+1$.
 Suppose that $u_1,\ldots,u_t$ are all vertices such that
$d(u_1,v_0)=\cdots=d(u_t,v_0)=2.$ Then $t\ge 1$. We may assume
without loss of generality that $v_1,\ldots,v_s$ are all vertices in
$\Gamma(v_0)$ adjacent to some vertex in the set
$\{u_1,\ldots,u_t\}$. Then $1\le s\le m$. Let
$V_1=\{v_0,v_1,\ldots,v_m\}, V_1'=V(G)-V_1$ and let $e(V_1V_1')$ be
the number of edges with one endpoint in $V_1$ and another endpoint
in $V_1'$. Since $G$ does not contain $T_n^{2}$, for $m\ge n-3$ each
$v_i(1\le i\le s)$ has one and only one adjacent vertex in the set
$\{u_1,\ldots,u_t\}$. Thus, for $m\ge n-3$ we must have
$e(V_1V_1')=s\ge t$.
\par If $m\ge n-1,$ since $G$
does not contain $T_n^2$ as a subgraph, we see that $d(v_i)\le 2$
for $i=1,\ldots,m$ and so
$e(G[V_1])=d(v_0)+e(G[v_{s+1},\ldots,v_m])\le m+\f{m-s}2$. Hence
$$\aligned e(G)&=e(G[V_1])+e(V_1V_1')+e(G-V_1)
\\&\le\f{3m-s}2+s+e(G-V_1)\le 2m+e(G-V_1).\endaligned$$
Suppose $m+1=k(n-1)+r$ with $k\in\Bbb N$ and $0\le r\le n-2$. Set
$G_1=kK_{n-1}\cup K_r$. Since $m+1\ge n$, by (2.2) we have $e(G_1)>
2(m+1)-1>2m$.
 Thus,
  $e(G_1\cup (G-V_1))=e(G_1)+e(G-V_1)>2m+e(G-V_1)\ge e(G)$. As
  $G_1$ does not contain
any copies of $T_n^2$ and
  $G$ is an extremal graph, we get a contradiction. Hence
$\Delta(G)=m\le n-2.$
\par Suppose $m=n-2$. As $G$ does not contain $T_n^2$ as a subgraph,
 we see that
$d(v_1)=\cdots=d(v_s)=2$ and so
 $e(G[V_1])\le
n-2+\b{n-2-s}2.$ Since $1\le s\le m=n-2\le 2n-8$, we have
$$\aligned e(G)&=e(G[V_1])+e(V_1V_1')+e(G-V_1)
\\&\le\b{n-2-s}2+n-2+s+e(G-V_1)
\\&=\f{(n-2)(n-1)-s(2n-7-s)}2+e(G-V_1)
\\&<\b{n-1}2+e(G-V_1)=e(K_{n-1}\cup (G-V_1)).\endaligned$$
 This is impossible since $G$ is an extremal graph.

\par By the above, $\Delta(G)\le n-3$.
We first assume $\Delta(G)=n-3.$ We claim that $d(v_i)\le n-4$ for
$i=1,2,\ldots,s$. If $i\in\{1,2,\ldots,s\}$ and $d(v_i)=n-3$, let
$u_j$ be the unique adjacent vertex of $v_i$ in $\{u_1,\ldots,u_t\}$
and let $V_2=\{v_0,v_1,\ldots,v_{n-3},u_j\}$. Then there is at most
one vertex adjacent to $u_j$ in $G-V_2$. Hence $e(G-V_1)\le
1+e(G-V_2)$. Since each $v_r\ (1\le r\le s)$ is adjacent to one and
only one vertex in $\{u_1,\ldots,u_t\}$ and $\Delta(G[V_1])\le n-3$,
we see that
$$e(G[V_1])=\f 12\sum_{r=0}^{n-3}d_{G[V_1]}(v_r)\le
\f{s(n-4)+(n-2-s)(n-3)}2 =\f{(n-2)(n-3)-s}2.$$ Note that $s\le
\Delta(G)=n-3$. From the above we deduce that
$$\align e(G)&=e(G[V_1])+e(V_1V_1')+e(G-V_1)
=e(G[V_1])+s+e(G-V_1)\\&\le e(G[V_1])+s+1+e(G-V_2) \le
\f{(n-2)(n-3)-s}2+s+1+e(G-V_2) \\& =\f{(n-2)(n-3)+s+2}2+e(G-V_2)\le
\f{(n-2)(n-3)+n-1}2+e(G-V_2)
\\&<\f{(n-1)(n-2)}2+e(G-V_2)=e(K_{n-1}\cup(G-V_2)).
\endalign$$
Since $K_{n-1}\cup(G-V_2)$ does not contain $T_n^2$ and $G$ is an
extremal graph, we get a contradiction. Hence the claim is true.
Thus, for $\Delta(G)=n-3$ we have $d_{G[V_1]}(v_i)\le n-5$ for
$i=1,2,\ldots,s$ and so
$$e(G[V_1])=\f
12\sum_{i=0}^{n-3}d_{G[V_1]}(v_i)\le
\f{s(n-5)+(n-2-s)(n-3)}2=\f{(n-2)(n-3)}2-s.\tag 3.2$$

\par Now we assume $p<2n-2$ and $p=n-1+r$. Then $1\le r<n-1.$
By the above, $\Delta(G)\le n-3$. Assume $\Delta(G)=n-3$. Then
 $|V(G-V_1)|=p-(n-2)=r+1<n$,
$\Delta(G-V_1)\le n-3$ and so $e(G-V_1)\le
\t{min}\{\b{r+1}2,\f{(r+1)(n-3)}2\}.$ Since $e(G[V_1])\le
\f{(n-2)(n-3)}2-s$ by (3.2), we deduce that
$$\aligned e(G)&=e(G[V_1])+e(V_1V_1')+e(G-V_1)
\\&\le \f{(n-2)(n-3)}2-s+s+\t{min}\Big\{\f{r(r+1)}2,\f{(r+1)(n-3)}2\Big\}
\\&=\cases \f{(n-2)(n-3)}2+\b{r+1}2&\t{if}\q
r\le n-3\\\f{(n-2)(n-3)}2+\f{(n-3)(n-1)}2&\t{if}\q r=n-2\endcases
\\&<\b{n-1}2+\b r2=e(K_{n-1}\cup K_r).\endaligned$$
This is impossible since $G$ is an extremal graph. Thus,
$\Delta(G)\le n-4$ for $p<2n-2.$ Now the proof is complete.

 \pro{Lemma 3.2} Let $p,n\in\Bbb N$, $p\ge n\ge 7$ and $G\in
\Ex(p;T_n^2)$. Suppose that $G$ is connected. Then $p<2n-2.$\endpro
Proof. By Lemma 3.1, we have $\Delta(G)\le n-3$ and so
$e(G)\le\f{(n-3)p}2.$ Assume that $p=k(n-1)+r$ with $k\in\Bbb N$ and
$r\in\{0,1,\ldots,n-2\}$. Let $G_1\in \Ex(n-1+r;K_{1,n-3})$. Then
$e(G_1)=[\f{(n-4)(n-1+r)}2]$ by Lemma 2.1. Hence, if
$(k-2)(n-1)-r\ge 2$, then

$$\aligned e((k-1)K_{n-1}\cup G_1)&=(k-1)\b{n-1}2+\Big[\f{(n-4)(n-1+r)}2\Big]
\\&=\f{(p-r-(n-1))(n-2)}2+\Big[\f{(n-4)(n-1+r)}2\Big]
\\&=\Big[\f{(n-3)p}2+\f{p-2r-2(n-1)}2\Big]
\\&=\Big[\f{(n-3)p}2+\f{(k-2)(n-1)-r}2\Big]
>\Big[\f{(n-3)p}2\Big]\ge e(G).\endaligned$$ This is impossible
since $(k-1)K_{n-1}\cup G_1$ does not contain $T_n^2$ as a subgraph
and $G\in \Ex(p;T_n^2)$. Thus $(k-2)(n-1)-r\le 1$. If $k=3$, then
$r=n-2$ and $p=3(n-1)+n-2=4n-5$ and so
$$\align e(G)&\le \Big[\f{(n-3)p}2\Big]\le \f{(n-3)(4n-5)}2
=\f{4n^2-17n+15}2
\\&<\f{4n^2-14n+12}2=3\b{n-1}2+\b{n-2}2
=e(3K_{n-1}\cup K_{n-2}).
\endalign$$
Since $3K_{n-1}\cup K_{n-2}$ does not contain $T_n^2$ and $G\in
\Ex(p;T_n^2)$, we get a contradiction. Thus
 $k\le 2$.
\par For $p=2(n-1)+r$ with $r\in\{0,1,2,n-4,n-3,n-2\}$ we see that
$r(n-2-r)<2n-2$ and so  $e(2K_{n-1}\cup K_r)=
\f{2(n-1)(n-2)+r(r-1)}2>\f{(n-3)(2n-2+r)}2\ge e(G)$.
 This contradicts the assumption $G\in \Ex(p;T_n^2)$. Now suppose
$p=2(n-1)+r$ with $3\le r\le n-5$. If $\Delta(G)\le n-4$, then
$e(G)\le \f{(n-4)p}2$.  From previous argument we have
$$\align e(K_{n-1}\cup G_1)&=\b{n-1}2
+\Big[\f{(n-4)(n-1+r)}2\Big]=\Big[\f{(n-3)p-r}2\Big]
\\&=\Big[\f{(n-4)p}2\Big]+n-1>\f{(n-4)p}2\ge e(G).\endalign$$
Since $K_{n-1}\cup G_1$ does not contain $T_n^2$ as a subgraph and
$G\in \Ex(p;T_n^2)$, we get a contradiction. Hence $\Delta(G)=n-3$.
Suppose $v_0\in V(G)$, $d(v_0)=n-3$,
$\Gamma(v_0)=\{v_1,\ldots,v_{n-3}\}$,
$V_1=\{v_0,v_1,\ldots,v_{n-3}\}$ and $V_1'=V(G)-V_1$.
      Suppose also that
        there are exactly $s$ vertices in $\Gamma(v_0)$ adjacent to some
    vertex in $V_1'$. Then $1\le s\le n-3$. By (3.2),
     $e(G[V_1])
     \le \f{(n-2)(n-3)}2-s$.
       As $G$ does not contain any copies of $T_n^2$,
    we see that $e(V_1V_1')=s$.
Since $|V(G-V_1)|=|V_1'|=p-(n-2)=n+r$ and $G-V_1$ does not contain
any copies of $T_n^2$ we see that $e(G-V_1)\le \ex(n+r;T_n^2)$.
\par
We claim that $$\ex(n+r;T_n^2)\le
\max\Big\{\frac{(n-4)(n+r)}2,\frac{(n-1)(n-2)+r(r+1)}2\Big\}$$ for
$3\le r\le n-5$.
 Let $G'\in \Ex(n+r;T_n^2)$. If $G'$ is connected, using
Lemma 3.1 we
  have $\Delta(G')\le n-4$ and so $e(G')\le\f{(n-4)(n+r)}2$. Now
  suppose that $G'$ is not connected. If
  $n_1,n_2\in\{1,2,\ldots,n-2\}$, from Lemma 2.3 we have
  $e(K_{n_1}\cup K_{n_2})<e(K_{n_1+n_2})$ for $n_1+n_2<n$ and
  $e(K_{n_1}\cup K_{n_2})<e(K_{n-1}\cup K_{n_1+n_2-(n-1)})$ for
  $n_1+n_2\ge n$. Thus, $G'=G_1'\cup G_2'$, where $G_1'$ and $G_2'$
  are components of $G'$ with $|V(G_1')|=p_1'<n-1$ and $|V(G_2')|=p_2'\ge
  n-1$. For $p_2'\ge n$ we have $p_1'\le r\le n-3$ and so $e(G_1')=
  \f{p_1'(p_1'-1)}2\le \f{(n-4)p_1'}2$. For $p_2'\ge n$ we also have
  $\Delta(G_2')\le n-4$ and so $e(G_2')\le \f{(n-4)p_2'}2$ by Lemma
  3.1. Hence for $p_2'\ge n$ we find that $e(G')=e(G_1')+e(G_2')\le
  \f{(n-4)p_1'}2+\f{(n-4)p_2'}2=\f{(n-4)(n+r)}2$. Now assume
  $p_2'=n-1$. Then $p_1'=r+1$ and $$e(G')=e(K_{n-1}\cup
  K_{r+1})=\f{(n-1)(n-2)+r(r+1)}2.$$
    Hence the claim is true and so
    $$e(G-V_1)\le \ex(n+r;T_n^2)\le \max\Big\{\frac{(n-4)(n+r)}2,
    \frac{(n-1)(n-2)+r(r+1)}2\Big\}.$$
   Thus,
$$\aligned e(G)&=e(G[V_1])+e(V_1V_1')+e(G-V_1)
\\&\le\frac{(n-2)(n-3)}2-s+s
+\max\Big\{\frac{(n-4)(n+r)}2,\frac{(n-1)(n-2)+r(r+1)}2\Big\}
\\&=\binom{n-1}2+
\max\Big\{\frac{(n-4)(n-1+r)-n}2,\frac{(n-1)(n-2)+r(r-1)}2-(n-2-r)\Big\}
\\&<\binom{n-1}2+
\max\Big\{\Big[\frac{(n-4)(n-1+r)}2\Big],\frac{(n-1)(n-2)+r(r-1)}2\Big\}
\\&=\max\Big\{e(K_{n-1}\cup G_1),e(2K_{n-1}\cup K_r)\Big\}.
\endaligned$$
This is impossible since $G$ is an extremal graph.
\par By the above we must have $k=1$ and so
$p=k(n-1)+r<2n-2$  as asserted.

\pro{Lemma 3.3} Let $p,n\in \Bbb N$, $p\ge n\ge 7$ and  $G\in
\Ex(p;T_n^{2})$. Suppose that $G$ is connected. Then $\Delta(G)=n-4$
and $e(G)=[\f{(n-4)p}2]$.
\endpro
Proof. By (3.1), $e(G)\ge [\f{(n-4)p}2]$. If $\Delta(G)\le n-5$,
using Euler's theorem we see that $e(G)=\f 12\sum_{v\in V(G)}d(v)\le
\f{(n-5)p}2.$ Hence $\f{(n-4)p-1}2\le[\f{(n-4)p}2]\le e(G)\le
\f{(n-5)p}2$. This is impossible. Thus $\Delta(G)\ge n-4$. By Lemmas
3.1 and 3.2, $\Delta(G)\le n-4$. Therefore $\Delta(G)=n-4$ and so
$e(G)=\f 12\sum_{v\in V(G)}d(v)\le\f{(n-4)p}2$. Recall that $e(G)\ge
[\f{(n-4)p}2]$. Then $e(G)=[\f{(n-4)p}2]$ as asserted.

\pro{Lemma 3.4} Let $p$ and $k$ be nonnegative integers, $p=5k+r$
and $r\in\{0,1,2,3,4\}$. Suppose that $G$ is a graph of order $p$
without $T_6^2$. Then $e(G)\le 2p-\f{r(5-r)}2$.
\endpro
Proof. Clearly $\Delta(T_6^2)=3$. We prove the lemma by induction on
$p$. For $p\le 5$ we have $e(G)\le \f{p(p-1)}2=2p-\f{r(5-r)}2$. Now
suppose that $p\ge 6$ and the lemma is true for all graphs of order
$p_0<p$ without $T_6^2$. If $\Delta(G)\le 3$, then $e(G)=\f
12\sum_{v\in V(G)}d(v)\le \f{3p}2\le 2p-3\le 2p-\f{r(5-r)}2$.
\par Suppose $\Delta(G)=m\ge 4$, $v_0\in V(G)$, $d(v_0)=m$,
$\Gamma(v_0)=\{v_1,\ldots,v_m\}$, $V_1=\{v_0,v_1,\ldots,v_m\}$ and
$V_1'=V(G)-V_1$. If $G[V_1]$ is a component of $G$, then
$e(G[V_1])=e(K_5)=10$ for $m=4$, and $e(G[V_1])\le m+\f m2=\f{3m}2$
for $m\ge 5$ since $d(v_i)\le 2$ for $i=1,2,\ldots,m$. By the
inductive hypothesis, $e(G[V_1'])\le 2(p-m-1)-\f{r_1(5-r_1)}2$,
where $r_1\in\{0,1,2,3,4\}$ is given by $p-m-1\e r_1\mod 5$. Thus,
for $m=4$ we have $e(G)=e(G[V_1])+e(G[V_1'])\le
10+2(p-5)-\f{r(5-r)}2=2p-\f {r(5-r)}2$, and for $m\ge 5$ we have
$e(G)=e(G[V_1])+e(G[V_1'])\le \f{3m}2+ 2(p-m-1)-\f{r_1(5-r_1)}2 \le
2p-2-\f m2\le 2p-3\le 2p-\f{r(5-r)}2$.
\par From now on we assume that $G[V_1]$ is not a component of $G$
and $m=\Delta(G)\ge 4$. Hence there is a vertex $u_1$ such that
$d(u_1,v_0)=2$ and $u_1v_1\in E(G)$ with no loss of generality. Then
$v_1v_i\notin E(G)$ for $i=2,3,\ldots,m$.
  For $m=4$ we see that $e(G[V_1])+e(V_1V_1')\le
4+4=8$. For $m\ge 5$ we see that $d(v_i)\le 2$ for $i=1,2,\ldots,m$
and so
 $e(G[V_1])+e(V_1V_1')\le \sum_{i=1}^md(v_i)\le 2m$. Hence,
 for $m\ge 4$ we have
$e(G)=e(G[V_1])+e(V_1V_1')+e(G[V_1'])\le 2m+e(G[V_1']).$ By the
inductive hypothesis, $e(G[V_1'])\le 2(p-m-1)-\f{r_1(5-r_1)}2$,
where $r_1\in\{0,1,2,3,4\}$ is given by $p-m-1\e r_1\mod 5$. Thus,
$e(G)\le 2m+2(p-m-1)-\f{r_1(5-r_1)}2=2p-2-\f{r_1(5-r_1)}2.$ For
$r_1\ge 1$ we have $e(G)\le 2p-2-2<2p-\f{r(5-r)}2$. For $r_1=0$ and
$r=0,1,4$ we have $e(G)\le 2p-2\le 2p-\f{r(5-r)}2$. Therefore, we
only need to consider the case $p\e m+1\e 2,3\mod 5$.
\par Now assume $p\e m+1\e 2,3\mod 5$ and
$\Gamma(u_1)-\{v_1,\ldots,v_m\}=\{w_1,\ldots,w_t\}$. As $m\ge 4$ we
have $m\ge 6$. Set $V_2=\{v_0,v_1,\ldots,v_m,u_1\}$ and
$V_2'=V(G)-V_2$. Since $d(v_i)\le 2$ for $i=1,2,\ldots,m$, we see
that
$$e(G)=e(G[V_2])+e(V_2V_2')+e(G[V_2'])
\le \sum_{i=1}^md(v_i)+t+e(G[V_2'])\le 2m+t+e(G[V_2']).$$ Note that
$p-m-2\e 4\mod 5$ and $e(G[V_2'])\le 2(p-m-2)-\f{4(5-4)}2$ by the
inductive hypothesis. We then have $e(G)\le 2m+t+2(p-m-2)-2=2p+t-6$.
For $t\le 3$ we get $e(G)\le 2p+t-6\le 2p-3=2p-\f{r(5-r)}2$. For
$t\ge 4$ set $V_3=\{v_0,v_1,\ldots,v_m,u_1,w_1, \ldots,w_t\}$
 and $V_3'=V(G)-V_3$. Since $d(v_i)\le 2$ for $i=1,2,\ldots,m$ and
 $d(w_j)\le 2$ for $j=1,2,\ldots,t$, using the inductive hypothesis
 we see that
 $$\align e(G)&=e(G[V_3])+e(V_3V_3')+e(G[V_3'])
 \le  \sum_{i=1}^md(v_i)+\sum_{j=1}^td(w_j)+e(G[V_3'])
 \\&\le 2m+2t+e(G[V_3'])
 \le 2m+2t+2(p-m-2-t)=2p-4\\&<2p-\f{r(5-r)}2.
 \endalign$$
\par By the above, the lemma has been proved by induction.

 \pro{Theorem 3.1}  Let $p,n\in\Bbb N,\ p\geq n-1\geq
4$ and $p=k(n-1)+r$, where $k\in \Bbb N$ and $r\in\{
0,1,\ldots,n-2\}$. Then
$$\aligned
\ex(p;T_n^2)&=
\max\Big\{\Big[\f{(n-2)p}2\Big]-(n-1+r),\f{(n-2)p-r(n-1-r)}2\Big\}
\\&= \cases [\f{(n-2)p}2]-(n-1+r) &\t{if}\ n\ge 16\ \t{and}\
3\le r\le n-6 \ \t{or if}\\&\ \ 13\le n\le 15 \ \t{and}\ 4\le r\le n-7,\\
\f{(n-2)p-r(n-1-r)}2&\t{otherwise.}\endcases\endaligned$$\endpro
Proof. Clearly $\ex(n-1;T_n^2)=e(K_{n-1})=\f{(n-2)(n-1)}2$. Thus the
result is true for $p=n-1$. Now we assume $p\ge n$. Since
$T_5^2\cong T_5'$, taking $n=5$ in [10, Theorem 3.1] we obtain the
result in the case $n=5$. For $n=6$ we see that $\ex(p;T_6^2)\ge
e(kK_5\cup K_r)=10k+\f{r(r-1)}2=2p-\f{r(5-r)}2.$ This together with
Lemma 3.4 gives the result in this case. Applying Lemmas 3.3, 2.3
and replacing $T_n^1$ with $T_n^2$ in the proof of Theorem 2.1 we
deduce the result for $n\ge 7$.

\pro{Corollary 3.1} Suppose $p,n\in\Bbb N$, $p\ge n\ge 5$ and
$n-1\nmid p$. Then $\f{(n-2)p}2-\f{(n-1)^2}8\le \ex(p;T_n^2)\le
\f{(n-2)(p-1)}2$.
\endpro

 \section* {4. The Ramsey number $r(T_n^i,T_n)$}
\pro{\q Lemma 4.1 ([9, Lemma 2.1])} Let $G_1$ and $G_2$ be two
graphs. Suppose $p\in\Bbb N, \  p\ge max\{|V(G_1)|$, $|V(G_2)|\}$
and $\ex(p;G_1)+\ex(p;G_2)<\b p2.$ Then $r(G_1,G_2)\le p.$\endpro
 Proof. Let $G$ be a graph of order $p$. If
$e(G)\le \ex(p;G_1)$ and $e(\overline{G})\le \ex(p;G_2)$, then
$\ex(p;G_1)+\ex(p;G_2)\ge e(G)+e(\overline{G})=\b p2.$ This
contradicts the assumption. Hence, either $e(G)>\ex(p;G_1)$ or
$e(\overline{G})>\ex(p;G_2).$ Therefore, $G$ contains a copy of
$G_1$ or $\overline{G}$ contains a copy of $G_2$. This shows that
$r(G_1,G_2)\le |V(G)|=p.$ So the lemma is proved.

 \pro{Lemma 4.2 ([9, Lemma 2.3])}
Let $G_1$ and $G_2$ be two graphs with $\Delta(G_1)=d_1\ge 2$ and
$\Delta(G_2)=d_2\ge 2$. Then
\par $(\t{\rm i})$ $r(G_1,G_2)\ge d_1+d_2-(1-(-1)^{(d_1-1)(d_2-1)})/2$.
\par $(\t{\rm ii})$ Suppose that $G_1$ is a connected graph of order
$m$ and $d_1<d_2\le m$. Then $r(G_1,G_2)\ge 2d_2-1\ge d_1+d_2$.
\par $(\t{\rm iii})$ If $G_1$ is a connected graph of order
$m$, $d_1\not=m-1$ and $d_2>m$, then $r(G_1,G_2)\ge d_1+d_2$.
\endpro

\pro{Theorem 4.1} Let $n\in\Bbb N$ and $i,j\in\{1,2\}$.
\par $(\t{\rm i})$ If $n$ is odd with $n\ge 17$, then
$r(T_n^i,T_n^j)=2n-7$.
\par $(\t{\rm ii})$ If $n$ is even with $n\ge 12$, then
$r(T_n^i,T_n^j)=2n-6$.
\endpro
Proof. Suppose $n\ge 12$. Since $\Delta(T_n^i)=\Delta(T_n^j)=n-3$,
from Lemma 4.2 we know that $r(T_n^i,T_n^j)\ge 2n-7$ for odd $n$,
and $r(T_n^i,T_n^j)\ge 2n-6$ for even $n$. If $n$ is odd with $n\ge
17$, using Theorems 2.1 and 3.1 (with $k=1$ and $r=n-6$) we see that
$$\align \ex(2n-7;T_n^i)&=\f{(n-2)(2n-7)-1}2-(2n-7)
<\f{(n-4)(2n-7)}2=\f 12\b{2n-7}2\endalign$$ and so
$\ex(2n-7;T_n^i)+\ex(2n-7;T_n^j)<\b{2n-7}2$. Thus, by Lemma 4.1 we
have $r(T_n^i,T_n^j)\le 2n-7$. Hence (i) is true. From Theorems 2.1
and 3.1 (with $k=1$ and $r=n-5$) we see that for $n\ge 12$,
$$\align \ex(2n-6;T_n^i)&=\f{(n-2)(2n-6)-4(n-5)}2=n^2-7n+16
\\&<n^2-\f {13}2n+\f{21}2=\f 12\b{2n-6}2\endalign$$
and so $\ex(2n-6;T_n^i)+\ex(2n-6;T_n^j)<\b{2n-6}2$. Thus, by Lemma
4.1 we have $r(T_n^i,T_n^j)\le 2n-6$. Hence $r(T_n^i,T_n^j)= 2n-6$
for even $n$, proving (ii).

 \pro{Lemma 4.3}
Let $n\in\Bbb N$,  $n\ge 5$ and $i\in\{1,2\}$. Let $G_n$ be a
connected graph of order $n$ such that $\ex(2n-5;G_n)<n^2-5n+4.$
Then $r(T_n^i,G_n)\le 2n-5.$\endpro
 Proof. By Theorems 2.1 and 3.1,
  $\ex(2n-5;T_n^i)=\f{(n-2)(2n-5)-3(n-4)}2=n^2-6n+11.$ Thus,
$$\ex(2n-5;G_n)+\ex(2n-5;T_n^i)<n^2-5n+4+n^2-6n+11
=\b{2n-5}2.$$ Appealing to Lemma 4.1 we obtain $r(T_n^i,G_n)\le
2n-5$.
 \pro{Lemma
4.4 ([10, Theorem 3.1])} Let $p,n\in\Bbb N$ with $p\geq n\geq 5$.
Let $r\in\{ 0,1,\ldots,n-2\}$ be given by $p\e r\mod{n-1}$. Then
$$\ex(p;T_n')= \cases \big[\f{(n-2)(p-1)-r-1}2\big] &\t{if}\
n\ge 7\ \t{and}\
2\le r\le n-4,\\
\f{(n-2)p-r(n-1-r)}2&\t{otherwise.}\endcases$$
\endpro

 \pro{Theorem 4.2} Let
$n\in\Bbb N$, $n\ge 8$ and $i\in\{1,2\}$. Then
$r(T_n^i,T_n')=r(T_n^i,T_n^*)=2n-5$.\endpro
 Proof. Let $T_n\in\{T_n',T_n^*)$. As
$2K_{n-3}$ does not contain any copies of $T_n^i$ and
$\overline{2K_{n-3}}=K_{n-3,n-3}$ does not contain any copies of
$T_n$, we see that $r(T_n^i,T_n)\ge 1+2(n-3)=2n-5$.
 Taking $p=2n-5$ and $r=n-4$ in Lemma 4.4 we
find that
$$\ex(2n-5;T_n')=\Big[\f{(n-2)(2n-6)-(n-4)-1}2\Big]
\le n^2-\f{11}2n+\f{15}2 <n^2-5n+4.$$ By [10, Theorem 4.1],
$$\ex(2n-5;T_n^*)=\f{(n-2)(2n-5)-3(n-4)}2=n^2-6n+11<n^2-5n+4.$$
Thus, applying Lemma 4.3 we obtain $r(T_n^i,T_n)\le 2n-5$. Hence
$r(T_n^i,T_n)=2n-5$ as asserted.
 \par{\bf Remark 4.1} Let $n\in\Bbb
N$, $n\ge 5$ and $i\in\{1,2\}$.
 From [5,
Theorem 3.1(ii)] we know that $r(K_{1,n-1},T_n^i)=2n-3$.
 \pro{Theorem 4.3} Let
$n\in\Bbb N$ and $i\in\{1,2\}$. Then $r(P_n,T_n^i)=2n-7$ for $n\ge
17$,  $r(P_{n-1},T_n^i)=2n-7$ for $n\ge 13$, $r(P_{n-2},T_n^i)=2n-7$
for $n\ge 11$ and $r(P_{n-3},T_n^i)=2n-7$ for $n\ge 8$.
\endpro
Proof. Suppose $n\ge 8$ and $s\in\{0,1,2,3\}$. From Lemma 4.2(ii) we
have $r(P_{n-s},T_n^i)\ge 2(n-3)-1=2n-7$. By (1.1),
$$\ex(2n-7;P_{n-s})=\cases
\f{(n-2)(2n-7)-5(n-6)}2=\f{(n-4)(2n-7)+16-n}2 &\t{if $s=0$,}
\\\f{(n-3)(2n-7)-3(n-5)}2=\f{(n-4)(2n-7)+8-n}2&\t{if $s=1$,}
\\\f{(n-4)(2n-7)-(n-4)}2&\t{if $s=2$,}
\\\f{(n-5)(2n-7)-(n-5)}2=\f{(n-4)(2n-7)+12-3n}2&\t{if $s=3$.}
\endcases$$
By Theorems 2.1 and 3.1,
$$\ex(2n-7;T_n^i)=\cases [\f{(n-4)(2n-7)}2]&\t{if $n\ge 16$,}
\\\f{(n-2)(2n-7)-5(n-6)}2=\f{(n-4)(2n-7)+16-n}2&\t{if $n<16$.}\endcases$$
\par For $n\ge 17,13,11$ or $8$ according as $s=0,1,2$ or $3$, from the above we
find $\ex(2n-7;P_{n-s})+\ex(2n-7;T_n^i)<\b{2n-7}2$ and so
$r(P_{n-s},T_n^i)\le 2n-7$ by Lemma 4.1. This completes the proof.

 \section*{5. The Ramsey number
$r(T_m^i,T_n)$ for $m<n$}

 \pro{Proposition 5.1
(Burr[1])} Let $m,n\in\Bbb N$ with $m\ge 3$ and $m-1\mid n-2.$ Let
$T_m$ be a tree on $m$ vertices. Then
$r(T_m,K_{1,n-1})=m+n-2.$\endpro
 \pro{Proposition 5.2 (Guo and Volkmann [5, Theorem 3.1])}
  Let $m,n\in\Bbb N, m\ge
3$ and $n=k(m-1)+b$ with $k\in\Bbb N$ and
$b\in\{0,1,\ldots,m-2\}\setminus \{2\}.$ Let $T_m\not=K_{1,m-1}$ be
a tree on $m$ vertices. Then $r(T_m,K_{1,n-1})\le m+n-3$. Moreover,
if $k\ge m-b$, then $r(T_m,K_{1,n-1})=m+n-3.$\endpro

\pro{Lemma 5.1 ([6, Theorem 8.3, pp.11-12])} Let $a,b,n\in\Bbb N$.
If $a$ is coprime to $b$ and $n\ge (a-1)(b-1)$, then there are two
nonnegative integers $x$ and $y$ such that $n=ax+by$.\endpro
 \pro{Theorem 5.1} Let $m,n\in\Bbb
N$, $n>m\ge 5$, $m-1\nmid n-2$
 and $i\in\{1,2\}$. Then $r(T_m^i,K_{1,n-1})=m+n-3$ or $m+n-4$.
Moreover, if $n\ge (m-3)^2+1$ or $m+n-4=(m-1)x+(m-2)y$ for some
nonnegative integers $x$ and $y$, then $r(T_m,K_{1,n-1})=m+n-3$ for
any tree $T_m\not=K_{1,m-1}$ of order $m$.
 \endpro
Proof. Let $T_m\not=K_{1,m-1}$ be a tree on $m$ vertices. From
Proposition 5.2 we know that $r(T_m,$ $K_{1,n-1})\le m+n-3$. By
Lemma 4.2(iii), $r(T_m^i,K_{1,n-1})\ge m-3+n-1$. Thus,
$r(T_m^i,K_{1,n-1})=m+n-3$ or $m+n-4$. If $n\ge (m-3)^2+1$, then
$m+n-4\ge (m-2)(m-3)$ and so $m+n-4=(m-1)x+(m-2)y$ for some
nonnegative integers $x$ and $y$ by Lemma 5.1. If
$m+n-4=(m-1)x+(m-2)y$ for $x,y\in\{0,1,2,\ldots\}$, setting
$G=xK_{m-1}\cup yK_{m-2}$ we see that $G$ does not contain any
copies of $T_m$ and $\overline G$ does not contain any copies of
$K_{1,n-1}$. Thus $r(T_m,K_{1,n-1})\ge 1+|V(G)|=m+n-3$. Now putting
 all the above together we obtain the theorem.
 \pro{Theorem 5.2} Let $m,n\in\Bbb
N$, $n>m\ge 6$,  $m-1\mid n-3$ and $i\in\{1,2\}$. Then
$r(T_m^i,T_n')=m+n-3$.
\endpro
Proof.  By Theorems 2.1 and 3.1,
$\ex(m+n-3;T_m^i)=\f{(m-2)(m+n-3)-(m-2)}2 <\f{(m-2)(m+n-3)}2.$ Thus
applying [9, Theorem 5.1] we obtain the conclusion.

\pro{Theorem 5.3} Suppose $i\in\{1,2\}$, $m,n\in\Bbb N$, $n>m\ge 7$
and $m-1\nmid (n-3)$. Then $m+n-5\le r(T_m^i,T_n')\le m+n-4$ and
$m+n-6\le r(T_m^i,T_n^*)\le m+n-4$. Moreover, if
$n=k(m-1)+b=q(m-2)+a$, $k,q\in\Bbb N$, $a\in\{0,1,\ldots,m-3\}$,
$b\in\{0,1,\ldots,m-2\}$ and  one of the following conditions holds:
$$\align &(\t{\rm 1})\q b\in\{1,2,4\},
\\ &(\t{\rm 2})\q b=0\qtq{and}k\ge 3,
\\&(\t{\rm 3})\q n\ge (m-3)^2+2,
\\&(\t{\rm 4})\q n\ge m^2-1-b(m-2),
\\ &(\t{\rm 5})\q\ a\ge 3\qtq{and}n\ge (a-4)(m-1)+4,\qq\qq\qq\qq\qq
\endalign$$
 then $r(T_m^i,T_n^*)=r(T_m^i,T_n')=m+n-4$.\endpro

Proof. By Lemma 4.2 we have $r(T_m^i,T_n')\ge m-3+n-2$ and
$r(T_m^i,T_n^*)\ge m-3+n-3$.   Since $m-1\nmid n-3$, we have
$m-1\nmid m+n-4$. From Corollaries 2.1 and 3.1 we find
$\ex(m+n-4;T_m^i)\le \f{(m-2)(m+n-5)}2$. Hence, by [9, Lemma 5.2] we
have $r(T_m^i,T_n')\le m+n-4$, and by [9, Lemma 4.2] we have
$r(T_m^i,T_n^*)\le m+n-4$. Now applying  [9, Theorems 4.4 and 5.4]
we deduce the remaining assertion.

\section*{6. The Ramsey number $r(G_m,T_n^j)$ for $m<n$}
\pro{\q Theorem 6.1} Let $m,n\in\Bbb N$, $m\ge 5$, $n\ge 8$, $n>m$
and $j\in\{1,2\}$. Then $r(K_{1,m-1},T_n^j)$ $=m+n-4$ or $m+n-5$.
Moreover, if $2\mid mn$, then $r(K_{1,m-1},T_n^j)=m+n-4$.\endpro
 Proof. From
Lemma 4.2 we deduce that $r(K_{1,m-1},T_n^j)\ge
m-1+n-3-(1-(-1)^{(m-2)(n-4)})/2$ $=m+n-4-(1-(-1)^{mn})/2$. So, it
suffices to prove that $r(K_{1,m-1},T_n^j)\le m+n-4$. By Lemma 2.1,
$\ex(m+n-4;K_{1,m-1})=[\f{(m-2)(m+n-4)}2]$. By Theorems 2.1 and 3.1,
we have $$\ex(m+n-4;T_n^j)=\Big[\f{(n-4)(m+n-4)}2\Big]\ \t{or}\
\f{(n-2)(m+n-4)-(m-3)(n-m+2)}2.$$ Since
$[\f{(m-2)(m+n-4)}2]+[\f{(n-4)(m+n-4)}2]
\le\f{(m+n-6)(m+n-4)}2<\b{m+n-4}2$ and
$$\align &\f{(m-2)(m+n-4)}2+\f{(n-2)(m+n-4)-(m-3)(n-m+2)}2
\\&=\f{(m+n-4)(m+n-5)-(m-4)(n-m-\f 2{m-4})}2<\b{m+n-4}2,\endalign$$
we see that $\ex(m+n-4;K_{1,m-1})+\ex(m+n-4;T_n^j) <\b{m+n-4}2$ and
so $r(K_{1,m-1},T_n^j)\le m+n-4$ by Lemma 4.1. This completes the
proof.

\pro{Theorem 6.2} Let $m,n\in\Bbb N$, $m\ge 4$, $n\ge 7$, $m-1\mid
n-4$
 and $j\in\{1,2\}$.
 \par $(\t{\rm i})$ If  $G_m$ is a connected graph of
order $m$ with $\ex(m+n-4;G_m)\le \f{(m-2)(m+n-5)}2$, then
$r(G_m,T_n^j)=m+n-4.$
\par $(\t{\rm ii})$  $r(T_m',T_n^j)=r(T_m^1,T_n^j)
=r(T_m^2,T_n^j)=m+n-4$ for $m\ge 5$, $r(T_m^*,T_n^j)=m+n-4$ for
$m\ge 6$,
 and
$r(P_m,T_n^j)=m+n-4$.
\endpro
Proof. Set $t=(n-4)/(m-1)$. Suppose that $G_m$ is a connected graph
of order $m$ with $\ex(m+n-4;G_m)\le \f{(m-2)(m+n-5)}2$. Then
clearly $\Delta(\overline{(t+1)K_{m-1}})=t(m-1)=n-4$. Thus,
$(t+1)K_{m-1}$ does not contain any copies of $G_m$ and
$\overline{(t+1)K_{m-1}}$ does not contain any copies of $T_n^j$.
Hence $r(G_m,T_n^j) \ge 1+ (t+1)(m-1)=m+n-4$. By Theorems 2.1 and
3.1,
$$\ex(m+n-4;T_n^j)=\Big[\f{(n-4)(m+n-4)}2\Big]\ \t{or}\
\f{(n-2)(m+n-4)-(m-3)(n-m+2)}2.$$ If
$\ex(m+n-4;T_n^j)=[\f{(n-4)(m+n-4)}2]$, then
$$\align &\ex(m+n-4;G_m)+\ex(m+n-4;T_n^j)\\&\le
 \f{(m-2)(m+n-5)+(n-4)(m+n-4)}2<\b{m+n-4}2.\endalign$$
If $\ex(m+n-4;T_n^j)=\f{(n-2)(m+n-4)-(m-3)(n-m+2)}2$, then
$$\align &\ex(m+n-4;G_m)+\ex(m+n-4;T_n^j)\\&\le
 \f{(m-2)(m+n-5)+(n-2)(m+n-4)-(m-3)(n-m+2)}2
 \\&=\b{m+n-4}2-\f{(m-4)(n-m+1)}2<\b{m+n-4}2.\endalign$$
Therefore, by Lemma 4.1 we always have $r(G_m,T_n^j)\le m+n-4$ and
hence $r(G_m,T_n^j)=m+n-4$. This proves (i).
\par Now consider (ii). Note that $m+n-4\e 1\mod{m-1}$.
By (1.1), we have
$\ex(m+n-4;P_m)=\f{(m-2)(m+n-5)}2$. By Lemma 4.4,
$\ex(m+n-4;T_m')=\f{(m-2)(m+n-5)}2$ for $m\ge 5$. By [10, Theorem
4.2], $\ex(m+n-4;T_m^*)=\f{(m-2)(m+n-5)}2$ for $m\ge 6$. By Theorems
2.1 and 3.1, $\ex(m+n-4;T_m^i)=\f{(m-2)(m+n-5)}2$ for $i\in\{1,2\}$
and $m\ge 5$.
 Thus from (i)
and the above we deduce (ii). The proof is complete.

\pro{Lemma 6.1} Let  $j\in\{1,2\}$,
$m,n\in\Bbb N$,
  $m\ge 7$ and $m-1\nmid n-4$.
   Assume $n=m+1\ge 12$ or $n\ge \t{max}\ \{m+2,19-m\}$.
\par $(\t{\rm i})$ If  $G_m$ is a connected graph of
order $m$ with $\ex(m+n-5;G_m)\le \f{(m-2)(m+n-6)}2$, then
$r(G_m,T_n^j)\le m+n-5.$
\par $(\t{\rm ii})$ For $T_m\in\{P_m,T_m',T_m^*,T_m^1,T_m^2\}$ we have
 $r(T_m,T_n^j)\le m+n-5$.
\endpro
Proof. Since $m+n-5=n-1+m-4$, by Theorems 2.1 and 3.1 we have
$$\align \ex(m+n-5;T_n^j)&=\Big[\f{(n-4)(m+n-5)}2\Big]
\\&\qtq{or}\f{(n-2)(m+n-5)-(m-4)(n-1-(m-4))}2.\endalign$$
If $n=m+1$, then $(m-4)(n-3-(m-4))=2(n-5)$. If $n\ge m+2$, then
$3\le m-4\le n-6$ and so
$(m-4)(n-3-(m-4))=(\f{n-3}2)^2-(m-4-\f{n-3}2)^2\ge
(\f{n-3}2)^2-(n-6-\f{n-3}2)^2=3(n-6)$. Thus,
$$\aligned &\f{(n-4)(m+n-5)+m-2}2-\f{(n-2)(m+n-5)-(m-4)(n-1-(m-4))}2
\\&=\f{(m-4)(n-3-(m-4))-2n+m}2
\\&\ge \cases\f{2(n-5)-2n+m}2=\f{m-10}2>0&\t{if $n=m+1\ge 12$,}
\\\f{3(n-6)-2n+m}2=\f{n-10+m-8}2>0&\t{if $n\ge \t{max}\ \{m+2,19-m\}$.}
\endcases\endaligned$$ Therefore, from the above we deduce that
$$\ex(m+n-5;T_n^j)<\f{(n-4)(m+n-5)+m-2}2.\tag 6.1$$
Hence, if $G_m$ is a connected graph of order $m$ with
$\ex(m+n-5;G_m)\le \f{(m-2)(m+n-6)}2$, then
$$\align &\ex(m+n-5;G_m)+\ex(m+n-5;T_n^j)\\&<
\f{(m-2)(m+n-6)}2+\f{(n-4)(m+n-5)+m-2}2 =\b{m+n-5}2.\endalign$$
Applying Lemma 4.1 we obtain (i).
\par Now we consider (ii).  Since $m-1\nmid (m+n-5)$, by
Corollaries 2.1 and 3.1 we have $\ex(m+n-5;T_m^i)\le
\f{(m-2)(m+n-6)}2$ for $i\in\{1,2\}$. By (1.1),
 $\ex(m+n-5;P_m)\le
\f{(m-2)(m+n-6)}2$. By Lemma 4.4, $\ex(m+n-5;T_m')\le
\f{(m-2)(m+n-6)}2$. By [10, Theorems 4.1-4.5], $\ex(m+n-5;T_m^*)\le
\f{(m-2)(m+n-6)}2$. Thus, from the above and (i) we deduce (ii).
This proves the lemma.

\pro{Theorem 6.3} Let $m\in\Bbb N$ and $j\in\{1,2\}$.
\par $(\t{\rm i})$ We have
$$r(T_m',T_{m+1}^j)=\cases 2m-4&\t{if $2\nmid m$ and $m\ge 9$,}
\\2m-5&\t{if $2\mid m$ and $m\ge 16$.}
\endcases$$
\par $(\t{\rm ii})$
If $n\in\Bbb N$, $m\ge 7$, $n\ge \t{max}\ \{m+2,19-m\}$ and
$m-1\nmid n-4$, then $r(T_m',T_n^j)=m+n-5$.
\endpro
Proof. We first assume $2\nmid m$ and $m\ge 9$. By Lemma 4.2(i), we
have $r(T_m',T_{m+1}^j) \ge m-2+m-2=2m-4$. By Lemma 4.4,
$\ex(2m-4;T_m')=\f{(m-2)(2m-4)-2(m-3)}2=m^2-5m+7$. By Theorems 2.1
and 3.1, $\ex(2m-4;T_{m+1}^j)=\f{(m-1)(2m-4)-4(m-4)}2=m^2-5m+10$.
Thus,
$$\align&\ex(2m-4;T_m')+\ex(2m-4;T_{m+1}^j)
\\&=m^2-5m+7+m^2-5m+10=2m^2-10m+17<2m^2-9m+10=\b{2m-4}2.
\endalign$$
Hence, by Lemma 4.1 we obtain $r(T_m',T_{m+1}^j)\le 2m-4$ and so
$r(T_m',T_{m+1}^j)=2m-4$.
 \par Now we assume $2\mid m$ and $m\ge 16$. By Lemma
4.2(i), $r(T_m',T_{m+1}^j) \ge m-2+m-2-1=2m-5$. By Lemma 4.4,
$\ex(2m-5;T_m')=[\f{(m-2)(2m-6)-(m-3)}2] =\f{2m^2-11m+14}2$. By
Theorems 2.1 and 3.1,
$\ex(2m-5;T_{m+1}^j)=[\f{(m-1)(2m-5)}2]-(2m-5)=\f{2m^2-11m+14}2$.
Thus,
$$ \ex(2m-5;T_m')+\ex(2m-5;T_{m+1}^j)
=2m^2-11m+14<2m^2-11m+15=\b{2m-5}2.$$ Hence, by Lemma 4.1 we obtain
$r(T_m',T_{m+1}^j)\le 2m-5$ and so $r(T_m',T_{m+1}^j)=2m-5$. This
proves (i).
\par Now we consider (ii). Suppose $n\in\Bbb N$, $m\ge 7$ and
$n\ge \t{max}\ \{m+2,19-m\}$. By Lemma 6.1(ii), $r(T_m',T_n^j)\le
m+n-5$. By Lemma 4.2, we have $r(T_m',T_n^j)\ge m-2+n-3$. Thus,
$r(T_m',T_n^j)=m+n-5$. This proves (ii). The proof is complete.

\pro{Theorem 6.4} Let $j\in\{1,2\}$, $m,n\in\Bbb N$, $m\ge 7$ and
$m-1\nmid n-4$. Suppose $n=m+1\ge 12$ or $n\ge \t{max}\
\{m+2,19-m\}$.
 Assume that $G_m\in\{P_m,T_m^*,T_m^1,T_m^2\}$ or $G_m$ is a connected graph of
order $m$ such that $\ex(m+n-5;G_m)\le \f{(m-2)(m+n-6)}2$.
 If $n\ge (m-3)^2+3$ or $m+n-6=(m-1)x+(m-2)y$ for
 some nonnegative integers $x$ and $y$, then
$r(G_m,T_n^j)=m+n-5.$
\endpro
Proof. If $n\ge (m-3)^2+3$, then $m+n-6\ge (m-2)(m-3)$ and so
$m+n-6=(m-1)x+(m-2)y$ for some $x,y\in\{0,1,2,\ldots\}$ by Lemma
5.1. Now suppose $m+n-6=(m-1)x+(m-2)y$, where
$x,y\in\{0,1,2,\ldots\}$.
 Set $G=xK_{m-1}\cup yK_{m-2}$. Then
$\Delta(\overline G)\le n-4$. Thus, $G$ does not contain any copies
of $G_m$ and $\overline G$ does not contain any copies of $T_n^j$.
Hence $r(G_m,T_n^j)\ge 1+|V(G)|=m+n-5$. On the other hand, by Lemma
6.1 we have $r(G_m,T_n^j)\le m+n-5$. Thus $r(G_m,T_n^j)=m+n-5$. This
proves the theorem.

 \pro{Corollary 6.1} Let $m,n\in\Bbb N$, $m\ge
7$, $m-1\mid n-b$, $b\in\{2,3,5\}$, $n\ge \t{max}\
 \{m+2,19-m\}$
 and $j\in\{1,2\}$.  Assume that $G_m\in\{P_m,T_m^*,T_m^1,T_m^2\}$
   or $G_m$ is a connected graph of
order $m$ with $\ex(m+n-5;G_m)\le \f{(m-2)(m+n-6)}2$. Then
$r(G_m,T_n^j)=m+n-5$.
\endpro
Proof. Set $k=(n-b)/(m-1)$. Then $k\in\Bbb N$.
 For $b=2$ we have $k\ge 2$.
Since
$$m+n-6=\cases (k-2)(m-1)+3(m-2)&\t{if $b=2$,}
\\(k-1)(m-1)+2(m-2)&\t{if $b=3$,}
\\(k+1)(m-1)&\t{if $b=5$,}\endcases$$
the result follows from  Theorem 6.4.

  \pro{Theorem 6.5} Let $m\in\Bbb N$, $m\ge
12$  and $i,j\in\{1,2\}$.
 Then
$$r(T_m^i,T_{m+1}^j)=r(T_m^*,T_{m+1}^j)=2m-5.$$
\endpro
Proof.  Let $T_m\in\{T_m^i,T_m^*\}$. By Theorems 2.1, 3.1 and [10,
Theorem 4.1],
$$\align& \ex(2m-5;T_m)=\f{(m-2)(2m-5)-3(m-4)}2,
\\& \ex(2m-5;T_{m+1}^j)= \f{(m-1)(2m-5)-5(m-5)}2\ \t{or}\
\Big[\f{(m-3)(2m-5)}2\Big].\endalign$$ Since
$\f{(m-2)(2m-5)-3(m-4)}2+\f{(m-3)(2m-5)}2
=\f{(2m-5)(2m-6)+7-m}2<\b{2m-5}2$ and
$$\align &\f{(m-2)(2m-5)-3(m-4)}2+\f{(m-1)(2m-5)-5(m-5)}2
\\&=2m^2-12m+26<2m^2-11m+15=\b{2m-5}2,
\endalign$$
we see that $\ex(2m-5;T_m)+\ex(2m-5;T_{m+1}^j)<\b{2m-5}2$.
 Hence, applying Lemma 4.1 we deduce that $r(T_m,T_{m+1}^j)\le
2m-5$. Since $\Delta(T_m)=m-3$ and $\Delta(T_{m+1}^j)=m-2$, by Lemma
4.2(i) we have $r(T_m,T_{m+1}^j)\ge m-3+m-2=2m-5$. Hence
$r(T_m,T_{m+1}^j)=2m-5$. This proves the theorem.
\par\q
\par{\bf Acknowledgements.}  The first author is supported by the
National Natural Science Foundation of China (grant No. 11371163),
and the second author is supported by the Fundamental Research Funds
for the Central Universities (grant No. 2014QNA58).

\end{document}